\documentclass[a4paper,11pt]{amsart}
\usepackage[latin1]{inputenc}
\usepackage{amsmath}
\usepackage{amsfonts}
\usepackage{amssymb}
\usepackage{amsthm}
\usepackage{mathrsfs}

\newcommand{\R}{{\mathbb{R}}}

\newcommand{\Z}{{\mathbb{Z}}}
\newcommand{\C}{\mathbb{C}}

\newcommand{\M}{\mathcal{M}}

\def\vecx{{\text{\boldmath$x$}}}
\def\vecy{{\text{\boldmath$y$}}}

\def\vecV{{\text{\boldmath$V$}}}

\def\vecm{{\text{\boldmath$m$}}}

\def\vecN{{\text{\boldmath$N$}}}

\def\vec0{{\text{\boldmath$0$}}}

\newcommand{\ve}{\varepsilon}

\newcommand{\sfrac}[2]{{\textstyle \frac {#1}{#2}}}

\newcommand{\SL}{\mathrm{SL}}

\newtheorem{prop}{Proposition}
\newtheorem{lem}{Lemma}
\newtheorem{thm}{Theorem}
\newtheorem{cor}{Corollary}

\newtheorem*{intro2}{Theorem 1$'$}

\newtheorem*{intro3}{Theorem 1$''$}

\theoremstyle{remark}
\newtheorem{remark}{Remark}

\newenvironment{Proof}{{\bf Proof. }}{\hfill$\square$\newline}

\begin{document}
\title[On the Poisson distribution of lengths of lattice vectors]{On the Poisson distribution of lengths of lattice vectors in a random lattice}
\author{Anders S\"odergren}
\address{Department of Mathematics, Uppsala University, Box 480,\newline
\rule[0ex]{0ex}{0ex}\hspace{8pt} SE-75106 Uppsala, Sweden\newline
\rule[0ex]{0ex}{0ex}\hspace{8pt} {\tt sodergren@math.uu.se}} 
\date{8 September 2010}

\maketitle

\begin{abstract}
We prove that the volumes determined by the lengths of the non-zero vectors $\pm\vecx$ in a random lattice $L$ of covolume $1$ define a stochastic process that, as the dimension $n$ tends to infinity, converges weakly to a Poisson process on the positive real line with intensity $\frac{1}{2}$. This generalizes earlier results by Rogers (\cite[Thm.\ 3]{rogers3}) and Schmidt (\cite[Satz 10]{schmidt}).
\end{abstract}



\section{Introduction}

Let $G_n=\SL(n,\R)$ and $\Gamma_n=\SL(n,\Z)$. We will be interested in the space $X_n=\Gamma_n\backslash G_n$ considered as the space of lattices of covolume $1$. We let $\mu_n$ denote the Haar measure on $G_n$, normalized so that it represents the unique right $G_n-$invariant probability measure on the homogeneous space $X_n$.

For a random lattice $L\in X_n$ we study the lengths of the pairs $\pm\vecx$ of non-zero vectors. Given a lattice $L\in X_n$, we order the non-zero vector lengths in $L$ as $0<\ell_1\leq \ell_2\leq \ell_3\leq\ldots$, where we count the common length of the vectors $\vecx$ and $-\vecx$ only once. For $j\geq1$, we define 
\begin{align*}
 \mathcal V_j:=\frac{\pi^{n/2}}{\Gamma(\frac{n}{2}+1)}\ell_j^n
\end{align*}
so that $\mathcal V_j$ is the volume of an $n$-dimensional ball of radius $\ell_j$. Finally, for $t\geq0$, we let 
\begin{align*}
 \tilde N_{t}(L):=\#\{j:\mathcal{V}_j\leq t\}.
\end{align*}

The first result concerning the statistics of the sequence $\{\ell_j\}_{j=1}^{\infty}$ is due to Rogers. In \cite{rogers3} he shows that, for fixed $t$, $\tilde N_{t}(L)$ has a distribution which converges weakly to the Poisson distribution with mean $\frac{t}{2}$ as $n\to\infty$. Clearly Rogers' result determines the limit distribution of each fixed $\ell_j$ as $n\to\infty$. In particular the volume $\mathcal{V}_1$ has an exponential distribution with expectation value $2$ as $n\to\infty$. In \cite{schmidt} Schmidt presents a more direct study of the distribution of $\ell_1$ using another method resulting in better error bounds as $n\to\infty$. 

There is a close connection between these results and lattice packings of spheres: For each fixed $n$ the random variable $\mathcal{V}_1$ has bounded range, i.e.\ $\sup_{L\in X_n}\mathcal V_1<\infty$. The quantity $\mathcal{D}_n:=2^{-n}\sup_{L\in X_n}\mathcal V_1$
is in fact the maximal density of a lattice sphere packing in $\R^n$. The search for upper bounds for the density of lattice (and also more general) sphere packings is a very well studied subject, yet optimal bounds are known only in very special cases. At present the best known bounds as $n\to\infty$ are due to Ball (\cite{ball}) and Kabatyanskii and Levenshtein (\cite{sphere}); these state that 
\begin{align*}
2(n-1)\zeta(n)2^{-n}\leq\mathcal{D}_n\leq2^{-n(0.599+o(1))}. 
\end{align*}
Note that the limit distribution result for $\mathcal{V}_1$ quoted above immediately implies that $2^n\mathcal{D}_n\to\infty$ as $n\to\infty$. This observation can be improved using knowledge on the rate of convergence to the limit distribution; in fact Schmidt can give a lower bound for $\mathcal{D}_n$ which differs from that of Ball only by a modest constant factor (cf.\ \cite[Thm.\ 11]{schmidt}).

The results of Rogers and Schmidt suggest very naturally that it should be possible to understand not only the distribution of each fixed $\ell_j$ but also the joint distribution of the entire sequence $\{\ell_j\}_{j=1}^{\infty}$ for a random lattice $L\in X_n$ as $n\to\infty$. The purpose of the present note is to point out that this is possible using Rogers' methods. Our main theorem states that, as $n\to\infty$, the volumes $\{\mathcal V_j\}_{j=1}^{\infty}$ behave like the points of a Poisson process on the positive real line with intensity $\frac{1}{2}$. 

\begin{thm}\label{mainthm}
For any fixed $N$, the $N$-dimensional random variable $(\mathcal V_1,\ldots,\mathcal V_N)$ converges in distribution to the distribution of the first $N$ points of a Poisson process on the positive real line with intensity $\frac{1}{2}$ as $n\to\infty$.  
\end{thm}

Or, put equivalently: 

\begin{intro2}
Let $\{N(t),t\geq0\}$ be a Poisson process on the positive real line with intensity $\frac{1}{2}$. Then the stochastic process $\{\tilde N_{t}(\cdot),t\geq0\}$ converges weakly to \mbox{$\{N(t),t\geq0\}$} as $n\to\infty$.
\end{intro2}

The author got interested in this problem during recent investigations of the value distribution of the Epstein zeta function. We recall that for $s>\frac{n}{2}$ and $L\in X_n$ the Epstein zeta function is defined by
\begin{align}\label{Epstein}
E(L,s)={\sum_{\vecm\in L}}'|\vecm|^{-2s},
\end{align}
where $'$ denotes that the zero vector should be omitted. $E(L,s)$ has an analytic continuation to $\C$ except for a simple pole at $s=\frac{n}{2}$. Now Theorem \ref{mainthm} together with equation \eqref{Epstein} suggest that for a fixed $c>\frac{1}{2}$ the distribution of the random variable $E(\cdot,cn)$ (suitably normalized) converges weakly to the distribution of a similar sum, defined in terms of the points of a Poisson process, as $n\to\infty$. This and related results will be presented in detail in \cite{jag}.

Let us also point out the close conceptual relation between our Theorem \ref{mainthm} and the main result (Theorem 1.1) by VanderKam in \cite{vander} (cf.\ also Sarnak \cite{sarnak}), which states that for a generic \textit{fixed} lattice $L\in X_n$, the local spacing statistics of the infinite sequence $\{\mathcal V_j\}_{j=1}^{\infty}$ exhibit Poissonian behavior (to a larger extent the larger the dimension). In precise terms:

\begin{thm}[VanderKam] 
For almost all $L\in X_n$ the $2$-, $3$-,..., and $[n/2]$-level correlations of the sequence $\{\mathcal V_j\}_{j=1}^{\infty}$ are Poissonian, i.e.\ for any given intervals $I_1,\ldots,I_{m-1}\subset\R$ \emph{($2\leq m\leq\frac{n}{2}$)} we have
\begin{align*}
\lim_{N\to\infty}\frac{\#\big\{(i_1,\ldots,i_m)\in(\Z^+)_*^m:\mathcal{V}_{i_1},\ldots,\mathcal{V}_{i_m}<N, \mathcal{V}_{i_1}-\mathcal{V}_{i_j}\in I_{j-1}
\big\}}{N}=2^{1-m}\prod_{j=1}^{m-1}|I_j|,
\end{align*}
where $(\Z^+)_*^m$ is the set of all $m$-tuples of distinct positive integers.
\end{thm} 

VanderKam's theorem is of particular interest in the context of the Berry-Tabor conjecture (cf.\ \cite[Cor.\ 1.2]{vander}). This conjecture states that for a ``generic'' completely integrable system the distribution of the local spacings between the eigenvalues of the Laplacian (suitably normalized) should match that of the points of a certain Poisson process on the positive real line. VanderKam's theorem gives evidence in support of this conjecture, since it implies that for almost all flat tori the distribution of the eigenvalues of the Laplacian agree, at least up to a certain number of correlations, with the distribution of the points of a Poisson process. In this context, let us note that our Theorem \ref{mainthm} can be viewed as a kind of ``low-eigenvalue-high-dimension'' analog of the Berry-Tabor conjecture for flat tori. Indeed, Theorem \ref{mainthm} can be reformulated as follows. (Note that the eigenvalues for the flat torus $\R^n/\Lambda$ are $4\pi^2|\ell|^2$ with $\ell$ belonging to the dual lattice $\Lambda^*$; ``desymmetrizing'' and renormalizing these to have mean spacing one we get the sequence $\{\frac{1}{2}\mathcal V_j\}_{j=1}^{\infty}$ for $\Lambda^*$.)

\begin{intro3}
For a random flat torus $\R^n/\Lambda$ with $\Lambda\in X_n$ the distribution of the 
non-zero eigenvalues (normalized to have mean-spacing one) converges to the distribution of the 
points of a Poisson process on the positive real line with intensity $1$ as $n\to\infty$.  
\end{intro3}

\section{Rogers' integration formula}

In this section we fix some notation concerning an integration formula that will be the major technical tool in the proof of our main result.

Let $1\leq k\leq n-1$ and let $\rho:(\R^n)^k\to\R$ be a non-negative Borel measurable function. In \cite{rogers1} Rogers considers the integral
\begin{align*}
\int_{X_n}\sum_{\vecm_1,\ldots,\vecm_k\in L}\rho(\vecm_1,\ldots,\vecm_k)\,d\mu_n(L),
\end{align*}
and shows that it equals a certain (positive) infinite linear combination of integrals of $\rho$ over various linear subspaces of $(\R^n)^k$. In this note we will be interested in the similar integral
\begin{align}\label{rogint}
\int_{X_n}\sum_{\vecm_1,\ldots,\vecm_k\in L\setminus\{\vec0\}}\rho(\vecm_1,\ldots,\vecm_k)\,d\mu_n(L).
\end{align}
It follows from Rogers' formula in \cite{rogers1} that the integral in \eqref{rogint} equals
\begin{align}\label{rogformula}
&\int_{\R^n}\cdots\int_{\R^n}\rho(\vecx_1,\ldots,\vecx_k)\,d\vecx_1\ldots d\vecx_k\\
&+\sum_{(\nu,\mu)}\sum_{q=1}^{\infty}\sum_{D}\Big(\frac{e_1}{q}\cdots\frac{e_m}{q}\Big)^n\int_{\R^n}\cdots\int_{\R^n}\rho\Big(\sum_{i=1}^m\frac{d_{i1}}{q}\vecx_i,\ldots,\sum_{i=1}^m\frac{d_{ik}}{q}\vecx_i\Big)\,d\vecx_1\ldots d\vecx_m.\nonumber
\end{align}
Here the outer sum is over all divisions $(\nu,\mu)=(\nu_1,\ldots,\nu_m;\mu_1,\ldots,\mu_{k-m})$ of the numbers $1,\ldots,k$ into two sequences $\nu_1,\ldots,\nu_m$ and $\mu_1,\ldots,\mu_{k-m}$ with $1\leq m\leq k-1$, satisfying
\begin{align}\label{division}
& 1\leq \nu_1<\nu_2<\ldots<\nu_m\leq k,\nonumber\\
& 1\leq\mu_1<\mu_2<\ldots<\mu_{k-m}\leq k,\\
& \nu_i\neq\mu_j, \text{ if $1\leq i\leq m$, $1\leq j \leq k-m$}.\nonumber
\end{align}
The inner sum in \eqref{rogformula} is over all $m\times k$ matrices $D$, with no column vanishing, with integer elements having greatest common divisor equal to $1$, and with
\begin{align*}
&d_{i\nu_j}=q\delta_{ij}, \hspace{10pt}i=1,\ldots,m,\,\,j=1,\ldots,m,\nonumber\\
&d_{i\mu_j}=0,\hspace{10pt}\text{ if }\,\mu_j<\nu_i,\,\,i=1,\ldots,m,\,\,j=1,\ldots,k-m.
\end{align*}
We call these matrices $(\nu,\mu)$-admissible. Finally $e_i=(\ve_i,q)$, $i=1,\ldots,m$, where $\ve_1,\ldots,\ve_m$ are the elementary divisors of the matrix $D$.


\begin{remark}
It follows from the conditions on the matrices $D$ above and \cite[Thm.\ 14.5.1]{hua} that in all cases we have $e_1=1$. In particular it follows that we always have
\begin{align*}
\Big(\frac{e_1}{q}\cdots\frac{e_m}{q}\Big)^n\leq q^{-n}. 
\end{align*}
\end{remark}

\section{The joint moments of $N_j(L)$}\label{Rogsection}

We will now consider the following situation: Let $k\geq1$ and fix $0<V_1\leq V_2\leq\ldots\leq V_k$. Given $n\geq k+1$ and $1\leq j\leq k$ we let $N_j(L)$ denote the number of non-zero lattice points of $L$ in the $n$-ball of volume $V_j$ centered at the origin. The goal of this section is to establish the following result concerning the joint moments of the functions $N_j(L)$:

\begin{thm}\label{poisson}
Let $k\geq1$ and fix $0<V_1\leq V_2\leq\ldots\leq V_k$. For each division $(\nu,\mu)$ in Rogers' formula let $M_{\nu,\mu}$ denote the number of $(\nu,\mu)$-admissible matrices $D$, with elements $d_{ij}\in\{0,\pm1\}$, having exactly one non-zero entry in each column. Then
\begin{align}\label{Thm2}
\mathbb E\Big(\prod_{j=1}^{k}N_j(L)\Big)&\to\prod_{j=1}^kV_j+\underset{(\nu,\mu)}{\sum} M_{\nu,\mu}\prod_{i=1}^mV_{\nu_i} 
\end{align}
as $n\to\infty$, where the sum is over all divisions possible in Rogers' formula.
\end{thm}

\begin{remark}
It is a straightforward calculation to verify that
\begin{align*}
 M_{\nu,\mu}=2^{k-m}\Big(\prod_{j=2}^{m-1}j^{\nu_{j+1}-\nu_j-1}\Big)m^{k-\nu_m}\,\,\,\,.
\end{align*}
Hence it is possible to rewrite the result in Theorem \ref{poisson} in a way not involving any reference to matrices. With notation as above we have
\begin{align*}
 \mathbb E\Big(\prod_{j=1}^{k}N_j(L)\Big)&\to\sum_{m=1}^{k}2^{k-m}\underset{1=\nu_1<\ldots<\nu_m\leq k}{\sum}\Big(\prod_{j=2}^{m-1}j^{\nu_{j+1}-\nu_j-1}\Big)m^{k-\nu_m}\prod_{i=1}^mV_{\nu_i}
\end{align*}
as $n\to\infty$.
\end{remark}

\subsection{The proof of Theorem \ref{poisson}}

We let $\rho_j$, $1\leq j\leq k$, be the characteristic function of the $n$-ball of volume $V_j$ centered at the origin. Then we can write 
\begin{align*}
 N_j(L)=\sum_{\substack{\vecm\in L\\\vecm\neq\vec0}}\rho_j(\vecm),\hspace{15pt} 1\leq j\leq k.
\end{align*}
Applying Rogers' formula yields
\begin{align}\label{pos1}
\mathbb E\Big(\prod_{j=1}^{k}N_j(L)\Big)&=\prod_{j=1}^k\int_{\R^n}\rho_j(\vecx)\,d\vecx+\underset{(\nu,\mu)}{\sum}\sum_{q=1}^{\infty}\sum_D\Big(\frac{e_1}{q}\cdots\frac{e_m}{q}\Big)^n\\
&\times\int_{\R^n}\cdots\int_{\R^n}\prod_{j=1}^k\rho_j\Big(\sum_{i=1}^m\frac{d_{ij}}{q}\vecx_i\Big)\,d\vecx_1\ldots d\vecx_m\,\,\, .\nonumber
\end{align}
It will turn out that the main contribution comes from the terms where $q=1$ and the matrix $D$ has entries $d_{ij}\in\{0,\pm1\}$ and exactly one non-zero entry in each column. We begin by determining the contribution to \eqref{pos1} from these terms.

\begin{prop}\label{Prop1}
If $1\leq k\leq n-1$, then 
\begin{align*}
\prod_{j=1}^k\int_{\R^n}\rho_j(\vecx)\,d\vecx+\underset{(\nu,\mu)}{\sum}\sum_D\int_{\R^n}\cdots\int_{\R^n}\prod_{j=1}^k\rho_j\Big(\sum_{i=1}^md_{ij}\vecx_i\Big)\,d\vecx_1\ldots d\vecx_m\\
=\prod_{j=1}^kV_j+\underset{(\nu,\mu)}{\sum}M_{\nu,\mu}\prod_{j=1}^mV_{\nu_j},  
\end{align*}
where the sum over $D$ is taken over all $(\nu,\mu)$-admissible matrices having entries $d_{ij}\in\{0,\pm1\}$ and exactly one non-zero entry in each column.
\end{prop}

\begin{Proof}
The matrices $D$ on this form satisfies by definition $d_{j,\nu_j}=1$, $1\leq j\leq m$. If we furthermore let $\lambda_{\ell}$ be such that $d_{\lambda_{\ell},\mu_{\ell}}=\pm1$, $1\leq \ell\leq k-m$, we find that
\begin{align*}
&\int_{\R^n}\cdots\int_{\R^n}\prod_{j=1}^k\rho_j\Big(\sum_{i=1}^md_{ij}\vecx_i\Big)\,d\vecx_1\ldots d\vecx_m\\
&=\int_{\R^n}\cdots\int_{\R^n}\prod_{j=1}^m\rho_{\nu_j}(\vecx_j)\prod_{\ell=1}^{k-m}\rho_{\mu_{\ell}}(\vecx_{\lambda_{\ell}})\,d\vecx_1\ldots d\vecx_m. 
\end{align*}
Since the volumes $V_j$ are increasing this integral equals
\begin{align*}
 \int_{\R^n}\cdots\int_{\R^n}\prod_{j=1}^m\rho_{\nu_j}(\vecx_j)\,d\vecx_1\ldots d\vecx_m=\prod_{j=1}^mV_{\nu_j},
\end{align*}
which gives the desired result. 
\end{Proof}

We next give a bound on the contribution from most of the remaining terms in \eqref{pos1}.

\begin{prop}\label{intermed3}
Let
\begin{align*}
 I(\nu,\mu,q,D)=\int_{\R^n}\cdots\int_{\R^n}\prod_{j=1}^k\rho_j\Big(\sum_{i=1}^m\frac{d_{ij}}{q}\vecx_i\Big)\,d\vecx_1\ldots d\vecx_m
\end{align*}
and let $M(D)$ be 
the largest value taken by any determinant of an $m\times m$-minor of $D$. Assume that $n>\underset{1\leq m\leq k-1}{\max}\big(m(k-m)+1\big)$. 
Then
\begin{align*}
\underset{(\nu,\mu)}{\sum}\sum_{q=1}^{\infty}\underset{\substack{D}}{\sum}\Big(\frac{e_1}{q}\cdots\frac{e_m}{q}\Big)^nI(\nu,\mu,q,D)=\underset{(\nu,\mu)}{\sum}\underset{\substack{D\\q=1\\M(D)=1}}{\sum}I(\nu,\mu,q,D)+R(k),   
\end{align*}
where the remainder term satisfies 
\begin{align*}
0\leq R(k)\ll 2^{-n}.
\end{align*}
The implied constant depends on $k$ and the volumes $V_1,V_2,\ldots,V_k$ but not on $n$. 
\end{prop}

\begin{Proof}
This is a trivial adaptation of \cite[Sec.\ 9]{rogers2} to our situation. (We bound $\rho_1,\rho_2,\ldots,\rho_{k-1}$ from above by $\rho_k$.) 
\end{Proof}


It is immediate to verify that matrices $D$ with $q=1$ and $M(D)=1$ have all entries $d_{ij}\in\{0,\pm1\}$. In particular all the matrices in Proposition \ref{Prop1} are on this form. In the remaining part of this section we will discuss the contribution to \eqref{pos1} coming from matrices $D$ with $q=1$ and $M(D)=1$ and at least one column containing more than one non-zero entry. Our starting point is the following lemma.

\begin{lem}\label{rogge}
For any $1\leq i< j<\ell\leq k$ and $\vecy\in\R^n$ we have
\begin{align*}
\int_{\R^n}\int_{\R^n}\rho_i(\vecx_1)\rho_j(\vecx_2)\rho_{\ell}(\pm\vecx_1\pm\vecx_2+\vecy)\,d\vecx_1d\vecx_2\leq2\big(\sfrac{3}{4}\big)^{\frac{n}{2}}V_{\ell}^2.
\end{align*}
\end{lem}

\begin{Proof}
Since
\begin{align*}
&\int_{\R^n}\int_{\R^n}\rho_i(\vecx_1)\rho_j(\vecx_2)\rho_{\ell}(\pm\vecx_1\pm\vecx_2+\vecy)\,d\vecx_1d\vecx_2\\
&\leq\int_{\R^n}\int_{\R^n}\rho_{\ell}(\vecx_1)\rho_{\ell}(\vecx_2)\rho_{\ell}(\pm\vecx_1\pm\vecx_2+\vecy)\,d\vecx_1d\vecx_2, 
\end{align*}
the lemma follows from \cite[Lemma 5]{rogers3}.
\end{Proof}

Following Rogers (\cite[Lemma 6]{rogers3}) we obtain the final estimate needed in the proof of Theorem \ref{poisson}.

\begin{lem}\label{lastrogge}
 Let $k\geq3$ and fix $m$ satisfying $2\leq m\leq k-1$. Let $(\nu,\mu)$ be a division of the numbers $1,\ldots,k$ satisfying \eqref{division} with our $m$. Let $D$ be a $(\nu,\mu)$-admissible matrix with $q=1$, $M(D)=1$ and at least one column containing more than one non-zero entry. Let $1\leq\ell\leq k-m$ be such that the leftmost column with more than one non-zero entry is $\mu_{\ell}$ and let $\lambda_1$ and $\lambda_2$ be minimal with the property that $\lambda_1<\lambda_2$ and $d_{\lambda_1,\mu_{\ell}}\neq0\neq d_{\lambda_2,\mu_{\ell}}$. Then
\begin{align}\label{last}
 \int_{\R^n}\cdots\int_{\R^n}\prod_{j=1}^k\rho_j\Big(\sum_{i=1}^md_{ij}\vecx_i\Big)\,d\vecx_1\ldots d\vecx_m\leq
2\big(\sfrac{3}{4}\big)^{\frac{n}{2}}V_{\mu_{\ell}}^2\prod_{j\in M}V_{\nu_j},
\end{align}
where $M=\{1,2,\ldots,m\}\setminus\{\lambda_1,\lambda_2\}$.
\end{lem}

\begin{Proof}
With the notation introduced above we get
\begin{align*}
&\int_{\R^n}\cdots\int_{\R^n}\prod_{j=1}^k\rho_j\Big(\sum_{i=1}^md_{ij}\vecx_i\Big)\,d\vecx_1\ldots d\vecx_m\\
&\leq\int_{\R^n}\cdots\int_{\R^n}\prod_{j\in M}\rho_{\nu_j}(\vecx_j)\bigg(\int_{\R^n}\int_{\R^n}\rho_{\nu_{\lambda_1}}(\vecx_{\lambda_1})\rho_{\nu_{\lambda_2}}(\vecx_{\lambda_2})\\
&\times\rho_{\mu_{\ell}}\Big(\pm\vecx_{\lambda_1}\pm\vecx_{\lambda_2}+\sum_{i=\lambda_2+1}^md_{i\mu_{\ell}}\vecx_i\Big)\,d\vecx_{\lambda_1}d\vecx_{\lambda_2}\bigg)\,\prod_{j\in M}d\vecx_j\\
&\leq2\big(\sfrac{3}{4}\big)^{\frac{n}{2}}V_{\mu_{\ell}}^2\prod_{j\in M}V_{\nu_j},
\end{align*}
where we have applied Lemma \ref{rogge} to the inner integral.
\end{Proof}

We note that for fixed $V_1,\ldots,V_k$ the right hand side of \eqref{last} is exponentially smaller (in $n$) than the right hand side of \eqref{Thm2}. We stress that also here the implied constant depends on $V_1,\ldots,V_k$. Since there are only finitely many matrices on the form in Lemma \ref{lastrogge} this, together with Propositions \ref{Prop1} and \ref{intermed3}, proves Theorem \ref{poisson}.

\section{The proof of Theorem \ref{mainthm}}\label{Latticepoisson}

Let us consider a Poisson process $\{N(t),t\geq0\}$ on the positive real line with constant intensity $\frac{1}{2}$. Thus $N(t)$ denotes the number of points falling in the interval $(0,t]$. We recall that $N(t)$ is Poisson distributed with expectation value $\frac{t}{2}$. We further let $T_1, T_2,T_3,\ldots$ denote the points of the process ordered in such a way that $0<T_1\leq T_2\leq T_3\leq\cdots$. 

\begin{prop}\label{svlem}
 Let $k\geq1$ and denote by $\mathcal{P}(k)$ the set of partitions of $\{1,\ldots,k\}$. For $1\leq j\leq k$ let $f_j:\R_{\geq0}\to\R$ be functions satisfying $\prod_{j\in B}f_j\in L^1(\R_{\geq0})$ for every nonempty subset $B\subseteq\{1,\ldots,k\}$. Then 
\begin{align}\label{almostking}
\mathbb E\Big(\prod_{j=1}^k\Big(\sum_{n=1}^{\infty}f_j(T_n)\Big)\Big)=\underset{P \in\mathcal{P}(k)}{\sum}2^{-\#P}\underset{B\in P}{\prod}\Big(\int_{0}^{\infty}\prod_{j\in B}f_j(x)\,dx\Big). 
\end{align}
\end{prop}

\begin{remark}
 When $k=1$ Proposition \ref{svlem} gives
\begin{align*}
 \mathbb E\Big(\sum_{n=1}^{\infty}f_1(T_n)\Big)=\frac{1}{2}\int_{0}^{\infty}f_1(x)\,dx,
\end{align*}
which agrees with \cite[p.\ 28 (3.18)]{king}.
\end{remark}

\begin{Proof}
It follows from Campbell's Theorem (\cite[p.\ 28]{king}; cf.\ also \cite[(3.14), (3.18) and Cor.\ 3.1]{king}) that
\begin{align}\label{campbell}
&\mathbb E\Big(\sum_{\substack{n_1,\ldots,n_v\in\Z^+\\n_i=n_j\Leftrightarrow i=j}}g_1(T_{n_1})g_2(T_{n_2})\cdots g_v(T_{n_v})\Big)\\
&=\prod_{j=1}^v\mathbb E\Big(\sum_{n=1}^{\infty}g_j(T_n)\Big)=2^{-v}\prod_{j=1}^v\int_{0}^{\infty}g_j(x)\,dx \nonumber
\end{align}
for all $v\in\Z^+$ and all $g_1,\ldots,g_v\in L^1(\R_{\geq0})$. Now the left hand side of \eqref{almostking} can be expressed as
\begin{align*}
 \underset{P \in\mathcal{P}(k)}{\sum}\mathbb E\Big(\sum_{\substack{n_1,\ldots,n_k\in\Z^+\\n_i=n_j\Leftrightarrow i\sim_Pj}}f_1(T_{n_1})f_2(T_{n_2})\cdots f_k(T_{n_k})\Big)
\end{align*}
where $\sim_P$ is the equivalence relation on $\{1,\ldots,k\}$ corresponding to $P$. Hence \eqref{almostking} follows from \eqref{campbell} applied for all $v=1,\ldots,k$ and with functions $f_B=\prod_{j\in B}f_j$, $B\in P$, in place of $g_1,\ldots,g_v$.
\end{Proof}

As in Section \ref{Rogsection} we let $k\geq1$ and fix $0<V_1\leq V_2\leq\ldots\leq V_k$. We will apply Proposition \ref{svlem} to the situation where $f_j=\chi_j$, $1\leq j\leq k$, where $\chi_j$ is the characteristic function of the interval $[0,V_j]$. In this case the proposition states that
\begin{align}\label{pos2}
\mathbb E\Big(\prod_{j=1}^kN(V_j)\Big)&=\mathbb E\Big(\prod_{j=1}^k\Big(\sum_{n=1}^{\infty}\chi_j(T_n)\Big)\Big)\nonumber\\
&=\underset{P \in\mathcal{P}(k)}{\sum}2^{-\#P}\underset{B\in P}{\prod}\Big(\int_{0}^{\infty}\prod_{j\in B}\chi_j(x)\,dx\Big)\\
&=\underset{P \in\mathcal{P}(k)}{\sum}2^{-\#P}\underset{B\in P}{\prod}V_{j_B}, \nonumber
\end{align}
where $j_B=\underset{j\in B}{\min}\,j$. 

\begin{remark}
The formula in \eqref{pos2} can also be obtained from the mixed cumulants of $N(V_1),\ldots,N(V_k)$ (see e.g.\ \cite[Prop.\ 6.16]{svante}). 
\end{remark}

We are interested in comparing the formula in \eqref{pos2} with the result of Theorem \ref{poisson}. We first need the following observation.

\begin{lem}\label{bijection}
Let $\mathcal{D}(k)$ be the set of matrices $D$ which have nonnegative entries and are on the form occurring in Theorem \ref{poisson}, together with the $k\times k$ indentity matrix $I_k$. When $D=I_k$ let $\nu_i=i$, $1\leq i\leq k$. Then there is a bijection $g:\mathcal{D}(k)\to\mathcal{P}(k)$ with the property that if $D\in\mathcal{D}(k)$ is an $m\times k$ matrix and $g(D)=P=\{B_1,\ldots,B_{\#P}\}$ then $\#P=m$ and $\{\nu_1,\ldots,\nu_m\}=\{j_{B_1},\ldots,j_{B_m}\}$.
\end{lem}

\begin{Proof}
For $D\in\mathcal{D}(k)$ of size $m\times k$ we set $B_i=\{j\mid d_{ij}\neq0\}$, $1\leq i\leq m$, and define $g(D)=\{B_1,\ldots,B_m\}$. It is clear that $g$ is a bijection with the desired properties. 
\end{Proof}

This lemma together with Theorem \ref{poisson} and the result in \eqref{pos2} immediately implies the following theorem.

\begin{thm}\label{poisson2}
Let $\tilde N_j(L)=\frac{1}{2}N_j(L)$ so that $\tilde N_j(L)$ denotes the number of pairs of non-zero lattice points of $L$ belonging to the $n$-ball centered at the origin having volume $V_j$, $1\leq j\leq k$. Then 
\begin{align*}
 \mathbb E\Big(\prod_{j=1}^{k}\tilde N_j(L)\Big)&\to \mathbb E\Big(\prod_{j=1}^kN(V_j)\Big)
\end{align*}
as $n\to\infty$.
\end{thm}

\begin{cor}\label{corconv}
Let $k\geq1$ and let $\vecV=(V_1,V_2,\ldots,V_k)$ for some fixed $0<V_1<V_2\ldots< V_k$. Consider the random vectors
\begin{align*}
\widetilde{\vecN}(L,\vecV)=\big(\tilde N_1(L),\ldots,\tilde N_k(L)\big) 
\end{align*}
and 
\begin{align*}
\vecN(\vecV)=\big(N(V_1),\ldots,N(V_k)\big). 
\end{align*}
Then $\widetilde{\vecN}(L,\vecV)$ converges in distribution to $\vecN(\vecV)$ as $n\to\infty$.
\end{cor}

\begin{Proof}
Since $N(V_j)$, $1\leq j\leq k$, is Poisson distributed with expectation value $\frac{V_j}{2}$ and the Poisson distribution is uniquely determined by its moments, it follows from \cite[Thm.\ 3]{petersen} that also the distribution of $\vecN(\vecV)$ is uniquely determined by its (joint) moments. The corollary now follows from this, Theorem \ref{poisson2} and a (in principle word by word) generalization to the present situation of \cite[Thm.\ 30.2]{billing}. 
\end{Proof}

\begin{remark}
When $k=1$ Corollary \ref{corconv} shows that $\tilde N_1(L)$ has a distribution which converges weakly to the Poisson distribution with mean $\frac{1}{2}V_1$ as $n\to\infty$. This was first proved by Rogers (\cite[Thm.\ 3]{rogers3}). 
\end{remark}

If we write $\tilde N_{V_1}(L)$ instead of $\tilde N_1(L)$ and then let $V_1=t$ we get a stochastic process $\{\tilde N_{t}(L),t\geq0\}$. It follows from Corollary \ref{corconv} that all finite dimensional distributions coming from this process converge to the corresponding finite dimensional distributions of the Poisson process $\{N(t),t\geq0\}$ as $n\to\infty$. Hence, by \cite[Thm.\ 12.6 and Thm.\ 16.7]{billconv}, the process $\{\tilde N_{t}(L),t\geq0\}$ converges weakly to the process $\{N(t),t\geq0\}$ as $n\to\infty$. This concludes the proof of Theorem \ref{mainthm}. 

\begin{remark}
The proof of Theorem \ref{mainthm} can be slightly simplified by considering factorial moments instead of ordinary moments (cf.\ \cite[Thm.\ 6.10]{svante}). However the present set-up will be useful also when dealing with the moments of the Epstein zeta function, cf.\ \cite{jag}. 
\end{remark}

\subsubsection*{Acknowledgement} The author is grateful to Andreas Strömbergsson for suggesting the problem and for many helpful discussions and to Svante Janson for valuable comments.

\end{document}